\documentclass[12pt]{amsart}

\batchmode

\begin{document}

\def\dbl{[\hskip -1pt[}
\def\dbr{]\hskip -1pt]}
\title {Finite type points on  subsets of  $\mathbb C^n$  }
\author{Ozcan Yazici}
\address{ Department  of Mathematics,  Middle East  Technical University,  06800 Ankara,  Turkey }
\email{oyazici@metu.edu.tr}




\def\Label#1{\label{#1}}

\def\1#1{\ov{#1}}
\def\2#1{\widetilde{#1}}
\def\6#1{\mathcal{#1}}
\def\4#1{\mathbb{#1}}
\def\3#1{\widehat{#1}}
\def\7#1{\mathscr{#1}}
\def\K{{\4K}}
\def\LL{{\4L}}

\def \MM{{\4M}}

\def \B{{\4B}^{2N'-1}}

\def \H{{\4H}^{2l-1}}

\def \F{{\4H}^{2N'-1}}

\def \LL{{\4L}}

\def\Re{{\sf Re}\,}
\def\Im{{\sf Im}\,}
\def\id{{\sf id}\,}

\def\s{s}
\def\k{\kappa}
\def\ov{\overline}
\def\span{\text{\rm span}}
\def\ad{\text{\rm ad }}
\def\tr{\text{\rm tr}}
\def\xo {{x_0}}
\def\Rk{\text{\rm Rk\,}}
\def\sg{\sigma}
\def \emxy{E_{(M,M')}(X,Y)}
\def \semxy{\scrE_{(M,M')}(X,Y)}
\def \jkxy {J^k(X,Y)}
\def \gkxy {G^k(X,Y)}
\def \exy {E(X,Y)}
\def \sexy{\scrE(X,Y)}
\def \hn {holomorphically nondegenerate}
\def\hyp{hypersurface}
\def\prt#1{{\partial \over\partial #1}}
\def\det{{\text{\rm det}}}
\def\wob{{w\over B(z)}}
\def\co{\chi_1}
\def\po{p_0}
\def\fb {\bar f}
\def\gb {\bar g}
\def\Fb {\ov F}
\def\Gb {\ov G}
\def\Hb {\ov H}
\def\zb {\bar z}
\def\wb {\bar w}
\def \qb {\bar Q}
\def \t {\tau}
\def\z{\chi}
\def\w{\tau}
\def\Z{{\mathbb Z}}
\def\phi{\varphi}
\def\eps{\epsilon}

\def \T {\theta}
\def \Th {\Theta}
\def \L {\Lambda}
\def\b {\beta}
\def\a {\alpha}
\def\o {\omega}
\def\l {\lambda}

\def \im{\text{\rm Im }}
\def \re{\text{\rm Re }}
\def \Char{\text{\rm Char }}
\def \supp{\text{\rm supp }}
\def \codim{\text{\rm codim }}
\def \Ht{\text{\rm ht }}
\def \Dt{\text{\rm dt }}
\def \hO{\widehat{\mathcal O}}
\def \cl{\text{\rm cl }}
\def \bS{\mathbb S}
\def \bK{\mathbb K}
\def \bD{\mathbb D}
\def \bC{\mathbb C}
\def \bL{\mathbb L}
\def \bZ{\mathbb Z}
\def \bN{\mathbb N}
\def \scrF{\mathcal F}
\def \scrK{\mathcal K}
\def \mc #1 {\mathcal {#1}}
\def \scrM{\mathcal M}
\def \cR{\mathcal R}
\def \scrJ{\mathcal J}
\def \scrA{\mathcal A}
\def \scrO{\mathcal O}
\def \scrV{\mathcal V}
\def \scrL{\mathcal L}
\def \scrE{\mathcal E}
\def \hol{\text{\rm hol}}
\def \aut{\text{\rm aut}}
\def \Aut{\text{\rm Aut}}
\def \J{\text{\rm Jac}}
\def\jet#1#2{J^{#1}_{#2}}
\def\gp#1{G^{#1}}
\def\gpo{\gp {2k_0}_0}
\def\emmp {\scrF(M,p;M',p')}
\def\rk{\text{\rm rk\,}}
\def\Orb{\text{\rm Orb\,}}
\def\Exp{\text{\rm Exp\,}}
\def\Span{\text{\rm span\,}}
\def\d{\partial}
\def\D{\3J}
\def\pr{{\rm pr}}

\def \CZZ {\C \dbl Z,\zeta \dbr}
\def \D{\text{\rm Der}\,}
\def \Rk{\text{\rm Rk}\,}
\def \CR{\text{\rm CR}}
\def \ima{\text{\rm im}\,}
\def \I {\mathcal I}

\def \M {\mathcal M}

\newtheorem{Thm}{Theorem}[section]
\newtheorem{Cor}[Thm]{Corollary}
\newtheorem{Pro}[Thm]{Proposition}
\newtheorem{Lem}[Thm]{Lemma}

\theoremstyle{definition}\newtheorem{Def}[Thm]{Definition}

\theoremstyle{remark}
\newtheorem{Rem}[Thm]{Remark}
\newtheorem{Exa}[Thm]{Example}
\newtheorem{Exs}[Thm]{Examples}

\numberwithin{equation}{section}

\def\bl{\begin{Lem}}
\def\el{\end{Lem}}
\def\bp{\begin{Pro}}
\def\ep{\end{Pro}}
\def\bt{\begin{Thm}}
\def\et{\end{Thm}}
\def\bc{\begin{Cor}}
\def\ec{\end{Cor}}
\def\bd{\begin{Def}}
\def\ed{\end{Def}}
\def\be{\begin{Exa}}
\def\ee{\end{Exa}}
\def\bpf{\begin{proof}}
\def\epf{\end{proof}}
\def\ben{\begin{enumerate}}
\def\een{\end{enumerate}}

\newcommand{\dbar}{\bar\partial}
\newcommand{\genmat}{\lambda}
\newcommand{\polynorm}[1]{{|| #1 ||}}
\newcommand{\vnorm}[1]{\left\|  #1 \right\|}
\newcommand{\asspol}[1]{{\mathbf{#1}}}
\newcommand{\Cn}{\mathbb{C}^n}
\newcommand{\Cd}{\mathbb{C}^d}
\newcommand{\Cm}{\mathbb{C}^m}
\newcommand{\C}{\mathbb{C}}
\newcommand{\CN}{\mathbb{C}^N}
\newcommand{\CNp}{\mathbb{C}^{N^\prime}}
\newcommand{\Rd}{\mathbb{R}^d}
\newcommand{\Rn}{\mathbb{R}^n}
\newcommand{\RN}{\mathbb{R}^N}
\newcommand{\R}{\mathbb{R}}
\newcommand{\bR}{\mathbb{R}}
\newcommand{\N}{\mathbb{N}}
\newcommand{\dop}[1]{\frac{\partial}{\partial #1}}
\newcommand{\vardop}[3]{\frac{\partial^{|#3|} #1}{\partial {#2}^{#3}}}
\newcommand{\br}[1]{\langle#1 \rangle}
\newcommand{\infnorm}[1]{{\left\| #1 \right\|}_{\infty}}
\newcommand{\onenorm}[1]{{\left\| #1 \right\|}_{1}}
\newcommand{\deltanorm}[1]{{\left\| #1 \right\|}_{\Delta}}
\newcommand{\omeganorm}[1]{{\left\| #1 \right\|}_{\Omega}}
\newcommand{\nequiv}{{\equiv \!\!\!\!\!\!  / \,\,}}
\newcommand{\bk}{\mathbf{K}}
\newcommand{\p}{\prime}
\newcommand{\tV}{\mathcal{V}}
\newcommand{\poly}{\mathcal{P}}
\newcommand{\ring}{\mathcal{A}}
\newcommand{\ringk}{\ring_k}
\newcommand{\ringktwo}{\mathcal{B}_\mu}
\newcommand{\germs}{\mathcal{O}}
\newcommand{\On}{\germs_n}
\newcommand{\mcl}{\mathcal{C}}
\newcommand{\formals}{\mathcal{F}}
\newcommand{\Fn}{\formals_n}
\newcommand{\autM}{{\Aut (M,0)}}
\newcommand{\autMp}{{\Aut (M,p)}}
\newcommand{\holmaps}{\mathcal{H}}
\newcommand{\biholmaps}{\mathcal{B}}
\newcommand{\autmaps}{\mathcal{A}(\CN,0)}
\newcommand{\jetsp}[2]{ G_{#1}^{#2} }
\newcommand{\njetsp}[2]{J_{#1}^{#2} }
\newcommand{\jetm}[2]{ j_{#1}^{#2} }
\newcommand{\glnc}{\mathsf{GL_n}(\C)}
\newcommand{\glmc}{\mathsf{GL_m}(\C)}
\newcommand{\glc}{\mathsf{GL_{(m+1)n}}(\C)}
\newcommand{\glk}{\mathsf{GL_{k}}(\C)}
\newcommand{\smC}{\mathcal{C}^{\infty}}
\newcommand{\anC}{\mathcal{C}^{\omega}}
\newcommand{\kC}{\mathcal{C}^{k}}


\begin{abstract} In \cite{D}, D'Angelo introduced the notion of points of finite type for a real hypersurface $M\subset \C^n$ and showed that the set of points of finite type  in $M$ is open. Later, Lamel-Mir  \cite{LM} considered a natural extension of D'Angelo's definition for an arbitrary set $M\subset \C^n$. Building on D'Angelo's work, we prove the openness of the set of points of finite type for any subset $M\subset \C^n$. 
 \end{abstract}

\maketitle

\section{Introduction}\Label{int} 

Let $M$ be a smooth real hypersurface of $\mathbb C^n$, $p\in M$ and $r$ be a defining function for $M$ such that $dr(p)$ does not vanish.  In his fundamental work \cite{D}, D'Angelo  introduced the notion of  type of $M$ at $p$ with respect to (possibly singular) holomorphic curves.  More precisely,  the type of $M$ at $p$ is defined by $$\Delta(M,p)=\sup_{\gamma\in\mathcal C} \frac{\nu(r\circ \gamma)}{\nu(\gamma)}$$ where  $\mathcal C$ is the set of non-constant holomorphic germs $\gamma$ at $0\in \mathbb C$ so that  $\gamma(0)=p$ and $\nu(r\circ\gamma)$ denotes the order of vanishing of the function $r\circ\gamma$ at $0.$ A point $p$ is called finite type if $\Delta(M,p)<\infty$.  In \cite{D}, D'Angelo proved the crucial property that the set of  points of finite type forms an open subset of the hypersurface $M$. This condition of finite type has been central in Catlin's work \cite{C} on subelliptic estimates for the $\bar \partial$-Neumann problem. See \cite{JM} for a more recent discussion of the relationship between finite type and subellipticity.

More recently, in their study of the ${ C}^\infty$ regularity problem for CR maps between smooth CR manifolds, Lamel and Mir \cite{LM} considered finite type points for arbitrary sets $M$ in $\C^n$. Building on D'Angelo's notion, they defined the type of an arbitrary set $M\subset \mathbb C^n$ as follows:
$$\Delta(M,p)=\sup_{\gamma\in\mathcal C} \inf_{r\in I_M(p)}\frac{\nu(r\circ \gamma)}{\nu(\gamma)}$$ where $I_M(p)$ is the set of germs at $p$ of real valued $C^{\infty}$ functions defined in a neighborhood of $p$ and vanishing on $M$ near $p$.  This definition of type coincides with the D'Angelo's definition for hypersurfaces, since the ideal $I_M(p)$ is generated by the defining function of $M$ when $M$ is a hypersurface.

As before, if $p\in M$, we say that $p$ is a point of finite type if $\Delta (M,p)<\infty$.  The goal of this  note is to show that D'Angelo's arguments can be used to establish the openness of the  set of points of finite type for any set $M\subset \mathbb C^n.$ 

Let  $\mathcal O_p$ denote the ring of holomorphic germs at $p$ in $\mathbb C^n.$ Following \cite{D}, for a proper ideal $I\subset\mathcal O_p$, we define \begin{eqnarray*}   \tau(I)=\sup_{\gamma\in\mathcal C}\inf_{\phi\in I}\frac{\nu(\phi\circ\gamma)}{\nu(\gamma)}\;\text{and}\; D(I)=\dim_{\mathbb C}(\mathcal O_p/I) \end{eqnarray*}  

If $I$ contains $q$ independent linear functions, then it follows from Theorem 2.7 in \cite{D} that $$\tau(I)\leq D(I)\leq (\tau(I))^{n-q}. $$

 For any  $r\in I_M(p)$ and $k\in\mathbb Z^+$, we denote by $r_k$ the Taylor polynomial of $r$ of order $k$ at $p$. We define $$\Delta(M_k,p)=\sup_{\gamma\in\mathcal C} \inf_{r\in I_M(p)}\frac{\nu(r_k\circ \gamma)}{\nu(\gamma)}. $$   Here $\nu(r_k\circ \gamma)$ denotes the multiplicity of  $r_k\circ \gamma(t)-r_k(p)$ at $0$.
By Proposition 3.1 in \cite{D}, we can decompose the polynomial $r_k$ as  $r_k=\Re (h^{r,k})+||f^{r,k} ||^2-||g^{r,k}||^2$ where  $h^{r,k}$ is a holomorphic polynomial and $f^{r,k}=(f_1^{r,k},\dots,f_N^{r,k})$ and  $g^{r,k}=(g_1^{r,k},\dots,g_N^{r,k})$ are holomorphic polynomial mappings. Here $N=N_k$ is independent of the polynomial and only depends on the degree $k$ and $n$. Let $\mathcal U(N_k) $ denote the group of unitary matrices on $\mathbb C^{N_k}$. For any $U\in \mathcal U(N_k) $, we denote   by $I(U,k,p)$, the ideal generated by the set $$\{h^{r,k}, \;\text{components of } \;f^{r,k}-Ug^{r,k}: r\in I_M(p) \}.$$  
We should note that the decomposition of $r_k$  is not unique and $I(U,k,p)$ depends on the choice of decomposition. Here, we use the one in the proof of Proposition 3.1 in \cite{D}.

\begin{Lem} \label{l1} Let $M\subset \mathbb C^n$ be a subset containing $p$.  If $\Delta (M_k,p)<k$ for some positive integer $k$, then $\Delta (M,p)=\Delta (M_k,p)$. 
\end{Lem}
\begin{proof} For any curve $\gamma\in \mathcal C$ and $r\in I_M(p)$ we write
 \begin{eqnarray}\label{1} r\circ\gamma =r_k\circ\gamma + e_{k}\circ \gamma \\ 
\label{2} \nu(e_k(\gamma))\geq (k+1)\nu(\gamma)
 \end{eqnarray} 
where $e_k=r-r_k$. Since  $\Delta (M_k,p)<k$,  for all $\gamma \in \mathcal C$ there exists $r^0\in I_M(p)$ such that $\frac{\nu(r_k^0\circ \gamma)}{\nu(\gamma)}<k$.  Then (\ref{1}) and (\ref{2}) imply that $\frac{\nu(r^0\circ\gamma)}{\nu(\gamma)}<k.$ By taking infimum over $r$ and supremum over $\gamma \in \mathcal C$ we obtain that $\Delta(M,p)\leq k.$ For all $\gamma \in \mathcal C$ and $\epsilon>0$  small enough, we can find $r^1\in I_M(p)$ such that $$\frac{\nu(r^1\circ \gamma)}{\nu(\gamma)}<  \inf_{r\in I_M(p)}\frac{\nu(r\circ \gamma)}{\nu( \gamma)}+\epsilon<k+1.$$ It follows from  (\ref{1}) and (\ref{2}) that $$\frac{\nu(r^1_k\circ \gamma)}{\nu( \gamma)}=\frac{\nu(r^1\circ \gamma) }{\nu( \gamma )}< \inf_{r\in I_M(p)}\frac{\nu(r\circ \gamma)}{\nu( \gamma)}+\epsilon.$$ As $\epsilon >0$ is arbitrary, $\displaystyle \inf_{r\in I_M(p) }\frac{\nu(r_k\circ \gamma)}{\nu( \gamma)}\leq \displaystyle \inf_{r\in I_M(p)}\frac{\nu(r\circ \gamma)}{\nu( \gamma)}. $ By taking supremum over $\gamma \in \mathcal C$, we get that $\Delta(M_k,p)\leq \Delta(M,p)$.
\\ 
We will show the other side of inequality in a similar way. Since we assume that $\Delta(M_k,p)<k$, for all $\gamma \in \mathcal C$ and $\epsilon >0$ there exists a $r^2\in I_M(p)$ such that   $$\frac{\nu(r_k^2\circ \gamma)}{\nu( \gamma)}<  \inf_{r\in I_M(p)}\frac{\nu(r_k\circ \gamma)}{\nu(\gamma)}+\epsilon<k.$$ It follows from   (\ref{1}) and (\ref{2})  that $$\frac{\nu(r^2\circ \gamma)}{\nu(\gamma)}=\frac{\nu(r^2_k\circ \gamma) }{\nu(\gamma) }< \inf_{r\in I_M(p)}\frac{\nu(r_k\circ \gamma)}{\nu(\gamma)}+\epsilon.$$ As $\epsilon >0$ is arbitrary, $\displaystyle \inf_{r\in I_M(p) }\frac{\nu(r\circ \gamma)}{\nu( \gamma)}\leq \displaystyle \inf_{r\in I_M(p)}\frac{\nu(r_k\circ \gamma)}{\nu( \gamma)}. $ By taking supremum over $\gamma \in \mathcal C$, we get that $\Delta(M,p)\leq \Delta(M_k,p)$.
\end{proof}

\begin{Lem}\label{l2}Let $M\subset \mathbb C^n$ be a subset containing $p$. Then    $$\displaystyle \sup_{U\in \mathcal U(N_k)} \tau(I(U,k,p))\leq \Delta(M_k,p)\leq 2 \sup_{U\in \mathcal U(N_k)} \tau(I(U,k,p)).$$
\end{Lem} 

\begin{proof} For any $r\in I_M(p)$, we can  write $r_k=\Re (h^{r,k})+||f^{r,k} ||^2-||g^{r,k}||^2$ where $h^{r,k}$ is a holomorphic polynomial function and $f^{r,k}$ $g^{r,k}$ are holomorphic polynomial mappings. We denote by $I(r,U, k,p)$ the ideal generated  by $h^{r,k}$ and the components of $f^{r,k}-Ug^{r,k}$ where $U\in \mathcal U(N_k)$. Since the supremum over $\gamma\in\mathcal C$ in the definition of $\tau(I(r,U,k,p))$ is attained for some $\gamma\in\mathcal C$ such that $h^{r,k}\circ \gamma=0$, it is enough to work with such curves.(See the proof of Theorem 1 on page 127 in \cite{D1}). As in the proof of Theorem 3.4 in \cite{JM}, $$\nu((f^{r,k}-Ug^{r,k})\circ\gamma)\leq \nu(r_k\circ\gamma)$$ for any $\gamma\in \mathcal C$ such that $h^{r,k}\circ\gamma=0$.  This implies that, $$\displaystyle \inf_{\phi\in I(r,U,k,p)} \nu(\phi\circ \gamma)\leq \nu(r_k\circ \gamma).$$ 
Dividing by $\nu(\gamma)$ and taking infimum over $r\in I_M(p)$, supremum over  $\gamma \in \mathcal C$ and over $U\in \mathcal U(N_k)$, we obtain that $$\displaystyle \sup_{U\in \mathcal U(N_k)} \tau(I(U,k,p))\leq \Delta(M_k,p).$$
 
 Let $r\in I_{M}(p)$, $\gamma\in \mathcal C$ and $2l+1\leq \nu(r_k\circ\gamma)\leq 2l+2$ for some $l\in \mathbb Z^+$. By Theorem 3.5 \cite{D}, there exists $U\in \mathcal U(N_k)$ such that $(h^{r,k}\circ \gamma)_{2l}=((f^{r,k}-Ug^{r,k})\circ\gamma)_l=0$ where $(h^{r,k}\circ \gamma)_{2l}$ and $((f^{r,k}-Ug^{r,k})\circ\gamma)_l$ are the Taylor polynomials of  $h^{r,k}\circ\gamma$ and  the components of $(f^{r,k}-Ug^{r,k})\circ\gamma$ of degree $2l$ and $l$, respectively. This implies that
 \begin{eqnarray} \label{m1}
 \frac{\nu(r_k\circ \gamma)}{\nu(\gamma)}\leq 2\inf_{\phi\in I(r,U,k,p)}\frac{\nu(\phi\circ \gamma)}{\nu(\gamma)}.\end{eqnarray}
When  $\nu(r_k\circ\gamma)=\infty$, by Theorem 3.5 in \cite{D},  there exists a $U\in\mathcal U(N_k)$ such that $\nu(\phi\circ\gamma)=\infty$ for all $\phi\in I(r,U,k,p)$ and the inequality (\ref{m1}) still holds. 
By taking infimum over $r\in I_M(p)$ and supremum over $\gamma\in \mathcal C$ and over $U\in \mathcal U(N_k)$, we obtain that $\Delta(M_k,p)\leq 2 \displaystyle \sup_{U\in \mathcal U(N_k)} \tau(I(U,k,p))$.
\end{proof} 

\begin{Thm}\label{main}  Let $M$ be  subset of $\mathbb C^n$ and $p_0$ be a point of finite type. Then there exists a neighborhood $V_0$ of $p_0$ such that \begin{eqnarray}\label{est} \Delta (M,p)\leq 2 (\Delta(M,p_0))^{n}, \end{eqnarray} for all $p\in V_0$. In particular, the set of points of finite type is an open subset of $M$. 
\end{Thm}

\begin{proof} We note that the coefficients of the generators of $I(U,k,p)$ depend smoothly on $p$. Then by Proposition 2.15 in \cite{D}, $ D(I(U,k,p))$  is an upper-semicontinuous function of $p$. Since $\mathcal U(N_k)$ is compact, upper-semicontinuity of $D$ implies that there exists a neighborhood $V_0$ of $p_0$ such that $$\sup_{U\in\mathcal U(N_k)}D(I(U,k,p))\leq \sup_{U\in\mathcal U(N_k)}D(I(U,k,p_0))  $$ for all $p\in V_0. $ By Theorem 2.7 in \cite{D}, \begin{eqnarray}\label{dt} \tau(I(U,k,p_0))\leq D(I(U,k,p_0))\leq (\tau(I(U,k,p_0)))^{n}. \end{eqnarray} 

By   Lemma  \ref{l2} we have $\Delta(M_k,p)\leq 2 \sup_{U\in\mathcal U(N_k)}\tau(I(U,k,p).$ Now we have the following chain of inequalities. For all $p\in V_0$, 

\begin{eqnarray*}  \Delta(M_k,p)&\leq &2\sup_{U\in\mathcal U(N_k)}\tau(I(U,k,p))\leq 2\sup_{U\in\mathcal U(N_k)}D(I(U,k,p)\\ \nonumber&\leq&  2\sup_{U\in\mathcal U(N_k)}D(I(U,k,p_0)  \leq 2\left( \sup_{U\in\mathcal U(N_k)}\tau(I(U,k,p_0))\right)^{n}\\&\leq&2 (\Delta(M_k,p_0))^{n}=2(\Delta(M,p_0))^{n}.
\end{eqnarray*}
The fifth inequality above follows from Lemma \ref{l2} and the  last equality follows from Lemma \ref{l1} for $k$ large enough. For any $k>2(\Delta(M,p_0))^{n}$, by Lemma \ref{l1},  $$\Delta(M,p)=\Delta(M_k,p)\leq  2(\Delta(M,p_0))^{n},$$ for all $p\in V_0.$
\end{proof}

\begin{Rem} If  $M$ is contained in a generic submanifold of real codimension $d$,  then we have a better estimate than  (\ref{est}).  In that case, there are $d$ local defining functions which vanish on $M$ near $p_0$ and their complex differentials are linearly independent  in a neighborhood $V_1$ of $p_0$.   After a change of coordinates  on $V_1$, we may assume that $I(U,k,p_0)$ contains $d$ independent  linear functions. Thus by Theorem 2.7 in \cite{D}, \begin{eqnarray}\label{dt} D(I(U,k,p_0))\leq (\tau(I(U,k,p_0)))^{n-d}. \end{eqnarray} 
We should note that after local biholomorphic change of coordinates, although the ideal $I(U,k,p_0)$ changes, the inequality (\ref{dt}) still holds.  Indeed, Corollary 3 on page 65 in  \cite{D1} implies that   $D(I)$ is invariant under a local  biholomorphic change of coordinates. Also, $\tau(I)$ is invariant under a local biholomorphic change of coordinates, (see the remarks after  Definition 8 on page 72 in \cite{D1}).Then it follows from a similar chain of inequalities as above that $ \Delta (M,p)\leq 2 (\Delta(M,p_0))^{n-d}$  for all $p\in V:=V_0\cap V_1.$
\end{Rem}

\begin{Rem} If a type  property fails for hypersurfaces, it also fails in higher codimension. For example, in \cite{D1}, D'Angelo gave an example of the hypersurface $$M=\{z\in\mathbb C^3: 2 \Re (z_3)+|z_1^2-z_2z_3|^2+|z_2|^4=0 \}$$ to show that the type $\Delta(M,p)$ is not an upper semicontinuous function of $p$. Here $\Delta(M,0)=4 $ and $\Delta(M,p)=8$ where $p=(0,0,ia).$   Following this example, we consider the set  $$M'=\{z\in \mathbb C^4: \Re(z_4)=   2 \Re (z_3)+|z_1^2-z_2z_3|^2+|z_2|^4=0  \}$$ of real codimension $2$ in $\mathbb C^4$ which has the same types as $M$. Hence, upper semicontinuity of type  fails  in higher codimension as well. 

\end{Rem}

\section{Points of finite q-type}
In  \cite{D}, D'Angelo defined the $q$-type of a hypersurface $M\subset \mathbb C^n$, which possibly contains $q-1$ dimensional complex analytic varieties. In a natural way, we can define $q$-type of an arbitrary subset $M\subset \mathbb C^n.$  The $q$-type of $M$ at $p\in M$ is defined by $$\Delta^q(M,p)=\inf\{\Delta(M\cap P,p): P\;\text{is any }\; n-q+1 \;\text{dimensional complex affine subspace of }\; \mathbb C^n  \}.  $$ More precisely, $$\Delta^q(M,p)=\inf\{ \Delta(\phi^{-1}(M),\phi^{-1}(p)): \phi: \mathbb C^{n-q+1}\to \mathbb C^n \;\text{is any linear imbedding} \}. $$

In \cite{C}, Catlin defined $q$-type, $D_q(M,p)$,  of a  hypersurface $M$ at $p$ by considering  generic $(n-q+1)$-dimensional complex affine  subspaces of $\mathbb C^n$ through $p$. When $q=1$, $D_1(M,p)=\Delta^1(M,p)=\Delta(M,p)$.  For a long time, Catlin's  $q$-type, $D_q(M,p)$, and D'Angelo's $q$-type,  $\Delta^q(M,p)$, were believed to be equal.  In \cite{F}, Fassina gave  examples of ideals and hypersurfaces to show that these two types are different when $q\geq2.$ This result also points out some mistakes in \cite{BN}, where the authors claimed that infimum in the definition of $\Delta^q(M,p)$ is achieved  by the  generic value with respect to choices of   $(n-q+1)$-dimensional complex affine subspaces of $\mathbb C^n$. Later, in \cite{BN1},  the authors corrected their results in \cite{BN}, by replacing $\Delta^q(M,p)$ with another invariant $\beta_q$ which is defined in terms of the generic value. 

With the D'Angelo's definition of $q$-type, openness of the set of points of  finite $q$-type easily follows from Theorem \ref{main}.
\begin{Cor} Let $M$ be a subset of $\mathbb C^n$. The set of points on $M$ for which $\Delta^q(M,p)<\infty$ is an open subset of $M.$  
\end{Cor}
\begin{proof}It follows from the definition of $q$-type that $$\{p\in M:\Delta^q(M,p)<\infty \}=  \cup_{\phi}A_{\phi}$$ where the union is taken over any imbedding $\phi:\mathbb C^{n-q+1} \to \mathbb C^n$ and $$A_{\phi}=\{z\in M:\Delta(\phi^{-1}(M),\phi^{-1}(p))<\infty\}.$$
Let $p_0\in A_{\phi}$.   By Theorem \ref{main}, there exists an open subset $U\subset \mathbb C^{n-q+1}$ around $\phi^{-1}(p_0)$ such that, for all $w\in U\cap \phi^{-1}(M)$, $\Delta(\phi^{-1}(M),w)<\infty.$ Since $\phi$ is an imbedding, $\phi(U)=V\cap\phi(\mathbb C^{n-q+1})$ for some open subset $V\subset \mathbb C^n$ containing $p_0$.  This implies that $M\cap V\subset A_{\phi}$ and $A_{\phi}$ is open in $M$. Thus $\{p\in M:\Delta^q(M,p)<\infty \}$ is an open set in $M$.  
\end{proof}

\noindent \textbf{Acknowledgments.} I am grateful to Prof. Nordine Mir for his suggestion to work on this question and for useful discussions on this subject. I would like to thank Prof. J.P. D'Angelo for his remarks and  suggestions which improved the  exposition of the paper.

\end{document}